\documentclass[12pt]{amsart}

\usepackage{geometry}
\usepackage{amssymb,latexsym, amscd, amsthm, amsmath}
\input xypic
\xyoption{all}

\newtheorem{theorem}{Theorem}[section]
\newtheorem{definition}[theorem]{Definition}
\newtheorem{proposition}[theorem]{Proposition}
\newtheorem{lemma}[theorem]{Lemma}

\newtheorem{remark}[theorem]{Remark}
\def\Proof{ \smallskip\noindent\rm {\sc Proof: }}
\def\Proofof{ \smallskip\noindent\rm {\sc Proof }of }

\begin{document}

\def\Tell{T_\ell}
\def\Zell{\mathbb Z_\ell}
\def\N{\mathbb{N}}
\def\Z{\mathbb{Z}}
\def\Q{\mathbb{Q}}
\def\Hzero{\operatorname{H}^0}
\def\Htate{{\widehat{\operatorname{H}}}^0}
\def\Hone{\operatorname{H}^1}
\def\Htwo{\operatorname{H}^2}
\def\Hmone{\widehat{\operatorname{H}}^{-1}}
\def\Ext{\operatorname{Ext}}
\def\Hom{\operatorname{Hom}}
\def\Ind{\operatorname{Ind}}
\def\Gal{\operatorname{Gal}}

\def\js{j_{!}}

\let\phi=\varphi
\let\cal=\mathcal

\parskip=2pt plus1pt

\title[Abstract $\ell$-adic $1$-motives and Tate's class] {{Abstract $\ell$--adic
$1$-motives and Tate's canonical class for number fields}}

  \author[Greither and Popescu]{Cornelius Greither and Cristian D. Popescu}

\address{Institut f\"ur Theoretische Informatik und Mathematik,
Universit\"at der Bundeswehr, M\"unchen,
85577 Neubiberg, Germany}
\email{cornelius.greither@unibw.de}

\address{Department of Mathematics, University of California, San Diego, La Jolla, CA 92093-0112, USA}
\email{cpopescu@math.ucsd.edu}

\keywords{Iwasawa theory;
Galois module structure; Tate sequences; \'Etale cohomology.}

\subjclass[2010]{11R23, 11R34, 11R33, 11R42, 11R70, 14G40}

\maketitle

\begin{abstract} In \cite{GP2} we constructed a new class of Iwasawa modules as $\ell$--adic realizations
of what we called {\it abstract $\ell$--adic $1$--motives} in the number field setting. We proved in loc. cit. that the new Iwasawa
modules satisfy an equivariant main conjecture. In this paper we link the new modules to the $\ell$--adified Tate canonical class,
defined by Tate in 1960 \cite{Ta1} and give an explicit construction of (the minus part of) $\ell$--adic Tate sequences for any Galois CM extension $K/k$
of an arbitrary totally real number field $k$. These explicit constructions are significant and useful in their own right but also due to their applications (via results
in \cite{GP2}) to a proof of the minus part of the far reaching Equivariant Tamagawa Number Conjecture for the Artin motive associated to the Galois extension $K/k$.

\end{abstract}

\section{Setup and preparation}\label{setup}

Let $K/k$ be a Galois extension of number fields of Galois group $G$.
Assume that $K$ is a CM field and that $k$ is totally real. We fix an odd prime $\ell$ and denote by $K_\infty$ and $k_\infty$ the cyclotomic
$\Bbb Z_\ell$--extensions of $K$ and $k$, respectively. We fix two finite, disjoint, $G$--invariant sets of primes
$S$ and $T$ in $K$, such that $S$ contains the ramification locus $S_{\rm ram}(K_\infty/k)$ of $K_\infty/k$ (in particular, it contains the set $S_\ell$ of all $\ell$--adic primes
and the set $S_\infty$ of all the archimedean primes) and $T$ contains at least two primes of distinct residual characteristics. We assume throughout that the {\it classical Iwasawa
$\mu$--invariant associated to $K_\infty$ and $\ell$ vanishes,} as conjectured by Iwasawa.

In earlier work \cite{GP2} we defined the category of ``abstract $\ell$--adic $1$--motives'' (which contains Deligne's category of Picard $1$--motives as a full subcategory)
and from the data $(K/k, S, T, \ell)$ as above we constructed a canonical abstract $\ell$--adic $1$--motive $\mathcal M:=\mathcal M_{S,T}^\ell(K/k)$. Its $\ell$--adic realization (Tate module)
$T_\ell(\mathcal M)$ which was defined in loc.cit. is a free $\Bbb Z_\ell$--module of finite rank which comes endowed with a natural $\Bbb Z_\ell[[\mathcal G]]$--module structure,
where $\mathcal G:={\rm Gal}(K_\infty/k)$. In fact, the unique complex conjugation automorphism $j$ of the CM field $K_\infty$ acts upon $T_\ell(\mathcal M)$ with eigenvalue $(-1)$,
so $T_\ell(\mathcal M)$ can be naturally viewed as a module over the quotient ring $\Bbb Z_\ell[[\mathcal G]]^-:= \Bbb Z_\ell[[\mathcal G]]/(1+j)$.
The main result in \cite{GP2} states the following.
\begin{theorem} Under the above hypotheses, the following hold.
\begin{enumerate}
\item ${\rm pd}_{\Bbb Z_\ell[[\mathcal G]]} T_\ell(\mathcal M)=1$.
\item If $G$ is abelian, then ${\rm Fit}_{\Bbb Z_\ell[[\mathcal G]]^-}T_\ell(\mathcal M)=\left(\Theta_{S,T}^\infty\right)$.
\end{enumerate}
\end{theorem}
Above, ``${\rm Fit}$'' denotes as usual the initial ($0$--th) Fitting ideal and $\Theta_{S,T}^\infty$  denotes a certain
equivariant $\ell$--adic $L$--function (a distinguished element of $\Bbb Z_\ell[[\mathcal G]]^-$) defined in loc.cit.
Part (2) of the above theorem is what we called an ``equivariant main conjecture'' and it is a $G$--equivariant refinement of the
classical Iwasawa Main Conjecture for arbitrary totally real number fields and odd primes $\ell$ proved by Wiles in \cite{Wi}. As shown in
\cite{GP2}, this refinement implies refined versions of the classical (imprimitive) Brumer-Stark and Coates-Sinnott conjectures.

From now on we will assume for simplicity that {\it the extensions $k_\infty/k$ and $K/k$ are linearly disjoint (over $k$).} This hypothesis
will be removed in Remark \ref{linear-disjointness-removed}. As a consequence of this hypothesis, Galois restriction induces a group isomorphism
$\mathcal G\simeq G\times \Gamma$, where $\Gamma:={\rm Gal}(K_\infty/K)\simeq {\rm Gal}(k_\infty/k)$. Consequently, we have ring isomorphisms
$\Bbb Z_\ell[[\mathcal G]]^-\simeq \Bbb Z_\ell[G]^-[[\Gamma]]\simeq \Lambda[G]^-$, where $\Lambda=\Bbb Z_\ell[[\Gamma]]$ is the usual Iwasawa algebra.
Consequently (see \cite{GP2} and the references therein) part (1) of the theorem above is equivalent to
$${\rm pd}_{\Bbb Z_\ell[G]}T_\ell(\mathcal M)=0,$$
i.e. $T_\ell(\mathcal M)$ is a finitely generated projective module over $\Bbb Z_\ell[G]$ (and over $\Bbb Z_\ell[G]^-$, obviously.)
As a consequence, if we fix a topological generator $\gamma$ of $\Gamma$, we obtain a perfect complex of $\Bbb Z_\ell[G]$--modules
$$   C^\bullet = [T_\ell(\mathcal M)  \stackrel{1-\gamma}{\longrightarrow}T_\ell(\mathcal M)], $$
concentrated in degrees $0$ and $1$. Of course,
the two cohomology groups of $C^\bullet$ are given by the $\Gamma$-invariants $T_\ell(\mathcal M)^\Gamma$
and $\Gamma$-coinvariants $T_\ell(\mathcal M)_\Gamma$ of $T_\ell(\mathcal M)$, respectively.
The goal of this paper
is to fully understand the two cohomology groups of $C^\bullet$, as well as the class of $C^\bullet$ in the
relevant $\Ext^2_{\Bbb Z_\ell[G]^-}(\ast, \ast)$. This goal will be stated much more precisely after the next remark.

\begin{remark}\label{remark-ffields} In \cite{GP1} we proved the exact analogue of Theorem 1.1 in the case where $K/k$ is a Galois extension
of global fields of characteristic $p>0$ (i.e. function fields) and $K_\infty$ (respectively $k_\infty$) is the maximal constant
field extension of $K$ (respectively $k$.) In that case there exists an actual geometric $1$--motive (Deligne's Picard $1$--motive)
$\mathcal M_{S,T}(K/k)$ whose $\ell$--primary part $\mathcal M_{S,T}(K/k)\otimes\Bbb Z_\ell$ gives the abstract $\ell$--adic
$1$--motive $\mathcal M_{S,T}^\ell(K/k)$, for all prime numbers $\ell$ (including $\ell=2, p$.) In that geometric context
there is no analogue of complex conjugation, so taking $(-1)$--eigenspaces does not make sense. Also, there is no analogue
of the sets $S_\ell$ or $S_\infty$ and, most importantly, the extension $K_\infty/K$ (which in that case is the maximal constant field extension of $K$)
is unramified.

Moreover, in \cite{GPff} we studied the function field analogue of the complex $C^\bullet$ and under a natural largeness
hypothesis on the set $S$ emerging from work of Tate (see below for details) we showed that there are $\Bbb Z_{\ell}[G]$--module isomorphisms
\begin{equation}\label{cohomology} {\rm H}^0(C^\bullet)\simeq U_{S,T}\otimes\Bbb Z_\ell, \qquad {\rm H}^1(C^\bullet)\simeq X_S\otimes\Bbb Z_\ell,\end{equation}
where $U_{S,T}$ is the group of $S$--units in $K$ which are congruent to $1$ modulo all primes in $T$ and
$X_S$ is the group of degree $0$ divisors in $K$ supported at $S$.

For all $G$--Galois extensions $K/k$ of global fields and data $(K/k, S, T)$ as above Tate \cite{Ta1, Ta2} defined a canonical
class $\tau_{K/k,S}\in \Ext^2_{\Bbb Z_\ell[G]}(X_S, U_S)$, for ``large'' $S$ and independent on $T$, where $U_S$ is the group of $S$--units in $K$ and $X_S$ is
as above. It turns out that under the ``largeness'' hypothesis (to be explained below),
the $\Bbb Z[G]$--module inclusion $\iota: U_{S,T}\to U_S$ induces a group isomorphism (more on this below)
\begin{equation}\label{iota-star}\iota_\ast: \Ext^2_{\Bbb Z[G]}(X_S, U_{S,T})\simeq \Ext^2_{\Bbb Z[G]}(X_S, U_{S}).\end{equation}
In the function field setting we proved in \cite{GPff} that if $c_{K/k, S, T}^{\ell}$ is the extension class of $C^\bullet$ and $\tau_{K/k, S}^{\ell}$ is the $\ell$-primary part of Tate's class, then
\begin{equation}\label{class}(\iota_\ast\otimes{\rm id}_{\Bbb Z_\ell}) (c_{K/k, S, T}^{\ell})=\tau_{K/k, S}^{\ell},\end{equation}
for all primes $\ell\ne p$. The same result should hold for $\ell=p$, but as explained
in loc. cit. we will address that case in a separate paper as the calculations would be somewhat different in nature, involving crystalline
rather than $\ell$--adic \'etale cohomology.
This way  we obtained in the function field setting a very explicit $\ell$--adic realization
$$  0 \to U_{S, T}\otimes\Bbb Z_\ell \to T_\ell(\mathcal M) \stackrel{1-\gamma}{\longrightarrow} T_{\ell}(\mathcal M) \to X_S\otimes\Bbb Z_{\ell} \to 0. $$
of a so--called Tate sequence (meaning that its middle terms are finitely generated, projective $\Bbb Z_\ell[G]$--modules and representing the $\ell$--adic Tate class via $\iota_\ast\otimes {\rm id}_{\Bbb Z_\ell}$.)

\end{remark}

The remark above makes it easier for us to state the goals of this paper more precisely: prove \eqref{cohomology} and \eqref{class} in the number field setting laid out above,
under the assumption that $\ell$ is an odd prime, with $U_{S,T}\otimes\Bbb Z_\ell$, $X_S\otimes\Bbb Z_\ell$ and $\tau_{K/k, S}^{\ell}$ replaced by
$(U_{S,T}\otimes\Bbb Z_\ell)^-$, $(X_S\otimes\Bbb Z_\ell)^-$ and $\tau_{K/k, S, T}^{\ell,\, -}$, respectively.
\medskip

As in \cite{GPff}, we will approach the question of linking $C^\bullet$
to the Tate class from two sides. On one hand we calculate
the $\Gamma$--invariants and $\Gamma$--coinvariants of $T_\ell(\cal M)$
directly, via Iwasawa theoretic methods, in sections \ref{Invariants}
and \ref{Coinvariants}. On the other hand, in section \ref{tate-link} (see Theorem \ref{mainresult}) we establish the desired
link between $C^\bullet$ and the Tate class $\tau_{K/k, S}^{\ell, -}$ via calculations in a certain
derived category and by relying in an essential way
upon deep results of Burns--Flach \cite{BF} and \cite{BuDoc}.
The reason why we insist on presenting the explicit calculations
of the cohomology of $C^\bullet$ is because the proof of Theorem \ref{mainresult}
relies on less explicit, not so easily transparent derived category arguments. It is satisfactory to see that the results obtained
via the two approaches agree at the cohomology level. We must admit that at present the explicit calculation in the
coinvariant case is somewhat laborious and not as smooth as the result one extracts
from the ``identification" with Tate's canonical class. However, it is definitely much more explicit.
\medskip

The rest of this section reviews additional notation and presents
some preparations.

For the construction of $\mathcal M$ and its $\ell$--adic realization $T_\ell(\mathcal M)$, the reader
should consult \cite{GP2}. For the definition of Tate's class $\tau_{K/k, S}$ the reader should consult \cite{Ta1, Ta2}.
In order to simplify notation we will let $\mathcal K:=K_\infty$.
For any algebraic field extension $N/K$, $S(N)$ denotes the set of places of $N$
above places in $S$, but often we will be sloppy in
context, just writing $S$ instead of $S(N)$.
 In the particular case $N=\cal K$ we write
$\cal S$ for $S(\cal K)$. The same notational convention
should be used for the set $T$, but for simplicity we will use $T$ for $T(N)$ and for $\mathcal T$ most of the time. No confusion will ensue.
The superscript minus always means the $(-1)$-eigenspace under
the unique complex conjugation of $\mathcal K$, as customary. As in Remark \ref{remark-ffields} above, $U_S$ denotes the group of $S$--units in $K$ and $U_{S,T}$
denotes its subgroup consisting of those $S$--units which are congruent to 1 modulo
every prime in $T$. For an algebraic extension $N/K$, $U_S(N)$ and $U_{S,T}(N)$ have similar meaning. If $X$ is a set and $O$ is a commutative ring,
then $O[X]$ will denote the free $O$--module of basis $X$. If $X$ happens to be a group (or a set endowed with an action by a group $H$), then $O[X]$ is viewed with its additional
group--ring structure (or $O[H]$--module structure). Note that since $K$ and $\mathcal K$ are CM, we have
$$(X_S\otimes\Bbb Z_\ell)^-=\Bbb Z_\ell[S]^-=\Bbb Z_\ell[S\setminus S_\infty]^-,\qquad  \Bbb Z_\ell[\cal S\setminus \cal S_\ell]^-=\Bbb Z_\ell[\cal S\setminus (\cal S_\ell\cup\cal S_\infty)]^-.$$

For an algebraic extension $N/K$, the group $cl_T(N)$ denotes the ray class group of $N$ with conductor
equal to the product of the prime ideals belonging to places in $T(N)$.
In less elaborate language, this is the group of all fractional ideals coprime to $T(N)$
modulo all principal ideals admitting a generator $u$ which
is congruent to 1 modulo all $v\in T(N)$. We let $cl(N)$ denote the usual class--group of $N$. For simplicity, we let $C^T(N):=(cl_T(N)\otimes\Bbb Z_\ell)^-$
and $C^T_\infty:=C^T(\mathcal K)$. We give similar meanings to $C(N)$ and $C_\infty$.

\begin{definition}\label{large-definition} The set $S$ is called large (respectively $\ell$--large) if $cl_T(K)$ (respectively $C^T(K)$)
is generated by ideal classes supported at primes in $S$.
\end{definition}

Note that Tate's definition of ``large'' involves the usual class--group $cl(K)$ instead of the ray--class group
$cl_T(K)$. However, the existence of a canonical surjective group morphism $cl_T(K)\twoheadrightarrow cl(K)$ shows that
``large'' in the sense of the definition above implies ``large'' in Tate's sense. Also, there is a well known canonical exact
sequence of $\Bbb Z[G]$--modules
\begin{equation}\label{basic-exact-seqeunce} 0 \to U_{S,T} \stackrel{\iota}{\to} U_S \to \kappa(T) \to cl_T(K)_S \to cl(K)_S\to 0\end{equation}
where $\kappa(T)=\oplus_{v\in T}\kappa(v)^\times$ (here $\kappa(v)$ is the residue field at $v$) and $cl_T(K)_S$ and  $cl(K)_S$
are the quotients of the corresponding ideal--class groups by the subgroups of $S$--ideal classes. It is well known (see \cite{GP2}, for example)
that ${\rm pd}_{\Bbb Z[G]}\kappa(T)=1$. Consequently, if $S$ is large then $cl_T(K)_S = cl(K)_S =0$ and $\iota$ induces the isomorphism
$\iota_\ast$ mentioned in \eqref{iota-star} above. Under the weaker ``$\ell$--largeness'' hypothesis this line of arguments yields the isomorphism
$(\iota_\ast\otimes{\rm id}_{\Bbb Z_\ell})^-$, which is in fact all that is needed for our goals.

We repeat our first goal: compute the modules $\Tell(\cal M)^\Gamma$ and
$\Tell(\cal M)_\Gamma$ directly.
The main problems we are going to encounter are caused
 by the set $S_\ell$ of $\ell$-adic places, which have no analog in the
function field case.
To guide us in our task, we recall from \cite{GP2} that there is
a canonical short exact sequence of $\Bbb Z_{\ell}[[\mathcal G]]$--modules
$$ 0 \to  \Tell(C_\infty^T) \to \Tell(\mathcal M) \to \Zell[\cal S\setminus\cal S_\ell]^- \to 0, $$
and we rely on the following
largely self-explanatory diagram arising from that s.e.s; the two dotted arrows indicate
the snake map. The resulting
6-term exact sequence of $\Gamma$-invariants and $\Gamma$-coinvariants,
connected by the snake map in the middle, is well visible
in this diagram and will be used later on. Here $\gamma$ is
a fixed generator of $\Gamma$.

$$ \xymatrix{
     & 0\ar[d] & 0\ar[d] & 0\ar[d]  \\
0 \ar[r]  & \Tell(C_\infty^T)^\Gamma \ar[r]\ar[d]  & \Tell(\mathcal M)^\Gamma
     \ar[r]\ar[d] & \Zell[\cal S\setminus\cal S_\ell]^{-,\Gamma} \ar@{.>}[r]\ar[d] &   \\
0 \ar[r] & \Tell(C_\infty^T) \ar[r]\ar[d]^{1-\gamma}
      & \Tell(\mathcal M) \ar[r]\ar[d]^{1-\gamma} &
     \Zell[\cal S\setminus\cal S_\ell]^- \ar[r]\ar[d]^{1-\gamma}  & 0 \\
0 \ar[r] & \Tell(C_\infty^T) \ar[r]\ar[d]  & \Tell(\mathcal M) \ar[r]\ar[d] &
     \Zell[\cal S\setminus\cal S_\ell]^- \ar[r]\ar[d]  & 0 \\
 \ar@{.>}[r] & \Tell(C_\infty^T)_\Gamma \ar[r]\ar[d]  & \Tell(\mathcal M)_\Gamma \ar[r]\ar[d] &
     \Zell[\cal S\setminus\cal S_\ell]^-_\Gamma \ar[r]\ar[d]  & 0 \\
     & 0 & 0 & 0            }$$

\section{Invariants}\label{Invariants}

We begin by dealing with the $\Gamma$-invariants. This is a relatively easy task in light of a very concrete interpretation
given to $T_{\ell}(\mathcal M)$ in \cite{GP2}, \S3. In this section and
the next, we make two blanket assumptions:
\begin{itemize}
  \item[(1)] $S$ is $\ell$--large, i.e. $C^T(K)$ is generated by the classes
of primes in $S$.
  \item[(2)] All primes in $S_\ell$ are totally
ramified in $K_\infty/K$.
\end{itemize}
The second assumption will be eliminated in section \ref{technical} below.

\begin{proposition}\label{invariants}
There is an isomorphism
$$ \phi_\infty:  \Tell(\mathcal M)^\Gamma \cong (\Zell \otimes_\Z U_{S,T})^-.$$

\end{proposition}

\Proof Recall from \S3 of \cite{GP2} that $\Tell(\cal M)\cong\varprojlim_\nu \cal M[\ell^\nu] $
and that there are canonical module isomorphisms
$$\cal M[\ell^\nu]\cong\left(\cal K^{(\ell^\nu)}_{\cal S,T}/\cal K_T^{\times \ell^\nu}\right)^-.$$
For simplicity, fix $\nu$, denote $m:=\ell^\nu$ and let $E_m:=\cal K^{(m)}_{\cal S,T}/\cal K_T^{\times m}$. Recall that
$$\cal K_T^{\times}:=\{x\in \cal K^\times\mid x\equiv 1{\,\rm mod\,} v, \forall\,v\in \mathcal T\}, \quad \cal K_{S, T}^{(m)}:=\{x\in \cal K_T^\times\mid {\rm div}_{\cal K}(x)=mD+D'\},$$
where ${\rm div}_{\cal K}(x)$ denotes the non-archimedean $\cal K$--divisor of $x$ and $D'$ is a divisor supported at $\mathcal S$. In plainer terms $\cal K_{S, T}^{(m)}$ consists of those elements
of $\mathcal K_T^\times$ whose divisors are multiples of $m$ away from $\mathcal S$.

(1) We claim
that  $E_m^\Gamma \cong (\cal K^{(m)}_{\cal S,T})^\Gamma/(\cal K_T^{\times m})^\Gamma$,
and that the denominator is simply $K_T^{\times m}$, where $K_T^\times$ is defined as above, but at the $K$--level.
Indeed, the second statement is clear
(raising to the power $m$ induces an isomorphism $\cal K_{\cal T}^\times\cong \cal K_T^{\times m}$, just as in loc.cit., since
there are no nontrivial $\ell$-power roots of unity in $\mathcal K_T^\times$, due to our assumptions on $T$).
For the first statement,
we need the vanishing of $\Hone(\Gamma, \cal K_T^{\times m})$. Again the exponent $m$
can be omitted, due to the isomorphism above. The vanishing follows, very similarly as in loc.cit., from Hilbert 90
and weak approximation. The ingredient which makes this work is the fact that
$T$ is unramified in the extension $\mathcal K/K$.

(2) By the previous step we have $E_m^\Gamma
\cong (\cal K^{(m)}_{\cal S,T})^\Gamma/K_T^{\times m}$. Now, we establish
a canonical isomorphism $$\pi_m:(\cal K^{(m)}_{\cal S,T})^\Gamma/K_T^{\times m}\cong U_{S,T}/U_{S,T}^{m}.$$

Take an element $x \in (\cal K^{(m)}_{\cal S,T})^\Gamma\subseteq K_T^\times$.
We have a unique writing ${\rm div}_{\cal K}(x) = m\cal D +\cal D'$ where $\cal D$ and $\cal D'$ are $\cal K$-divisors with
$\cal D'$ supported on $\cal S$ and $\cal D$ supported away from $\cal S$. Since $\cal K/K$ is unramified away from $S$ and $x\in K_T^\times$,
we also have ${\rm div}_K(x)=mD+D'$ with $K$--divisors $D$ and $D'$ supported away from and on S, respectively.
Using the first of our blanket hypotheses we get that $D = {\rm div}_K(y)+D''$ with $y\in K_T^\times$
and $D''$ supported on $S$. Hence ${\rm div}_K(xy^{-m})=mD''+D'$ is supported on $S$, and
therefore $xy^{-m} \in U_{S,T}$. We let $\tilde \pi_m(x):=\widehat{xy^{-m}}$. It is
easy to see that $\tilde \pi_m:  (\cal K^{(m)}_{\cal S,T})^\Gamma \to U_{S,T}/U_{S,T}^m$ is
well defined and onto, and also easily checked that the kernel is exactly
$K_T^{\times m}$. Therefore it induces the desired isomorphism $\pi_m$.

(3) After a compatibility check for the $\pi_m$'s
and passing to the projective limit, $\pi_\infty^- =\varprojlim_\nu \pi_{\ell^\nu}^-$
gives the desired isomorphism
$$\Tell(\cal M)\cong \varprojlim_\nu E_{\ell^\nu}^- \cong \left (\Zell \otimes_\Z U_{S,T}\right)^-.$$ We leave these details
to the interested reader.\ Q.E.D.


\section{Coinvariants}\label{Coinvariants}

Now, we turn to the calculation of $\Gamma$-coinvariants of $T_\ell(\cal M)$. We remind the reader that the
 assumptions (1) and (2), see beginning
of Section \ref{Invariants}, are in force.
The desired isomorphism $T_\ell(\mathcal M)_{\Gamma}\cong \Zell[S]^-$ will result via a simple homological
algebra lemma (Lemma \ref{perm-modules}) from Thm.~\ref{coinvses} (ii) below which yields
a short exact sequence
$$  0 \to \Zell[S_\ell]^-  \to \Tell(\mathcal M)_\Gamma
   \to \Zell[S\setminus S_\ell]^- \to 0. $$
Unfortunately there does not seem to be a simple
proof of the existence of this sequence. We begin with some notation and some fairly easy
auxiliary results. Then we present the calculation of the
coinvariants modulo three lemmas (one of which is highly technical),
and finally we proceed to prove the lemmas.

Let $K_n$ be the unique intermediate field of $\cal K/K$ with $[K_n:K]=\ell^n$. Let $\Gamma_n={\rm Gal}(K_n/K)$ and let $\gamma_n \in \Gamma_n$ be the image of
the generator $\gamma$ of $\Gamma$ via Galois restriction. Let $d$ be the $\Zell$-rank
of $\Zell[S_\ell]^-$ (note that this is unchanged
if $S_\ell$ is replaced by $S_\ell(K_n)$ or $\cal S_\ell$ due to our blanket hypotheses).

We remind the reader that $C^T(N)$ (respectively $C(N)$) is shorthand for the minus part
of the $\ell$-part of the ray class group $cl_T(N)$ (respectively class group $cl(N)$), for any
appropriate field $N$. (Usually $N$ is one of the fields $K_n$.) It is well known (see \cite{GP2}, for example) that
the canonical maps $C^T(K_n)\to C^T(K_{n+1})$ and $C(K_n)\to C(K_{n+1})$ are injective and that $C^T_\infty=\bigcup_n C_T(K_n)$ and
$C_\infty=\bigcup_n C(K_n)$.

Let $D^T(N) \subset C^T(N)$ be the subgroup generated
by the classes of the prime ideals in $N$ dividing $\ell$. It is easy to see that
$$   \mid {\rm Im}\Bigl(D^T(N) \to\, C^T(K_n)/C^T(K)\Bigr)\mid\,  \le  \ell^{nd}. $$
Note that it is legitimate to consider $C^T(K)$ as a subgroup of $C^T(K_n)$.

\begin{lemma}\label{one}
The preceding inequality is an equality, that is:
$$   \mid {\rm Im}\Bigl(D^T(N) \to \,C^T(K_n)/C^T(K)\Bigr)\mid\,   =  \ell^{nd}, \quad\text{for all $n$.} $$
\end{lemma}

\Proof Let $b_1,\ldots,b_d$ be a $\Zell$--basis of $\Zell[S_\ell(K_n)]^-$ where each $b_i$ has
the form $(1-j)\mathfrak p$ for some prime $\mathfrak p|\ell$ in $K_n$. (The letter $j$ means
complex conjugation of course; we have to take
exactly those $\mathfrak p$ that split from $K^+$ to $K$.) There is a map
$$ \phi_n^T : (\Z/\ell^n)^d \to  C^T(K_n)/C^T(K) $$
sending the $i$-th basis vector $e_i$ of the left-hand module
to the class of $b_i$. It is well-defined since $\ell^nb_i$ comes
from an ideal of $K$. The image of $\phi_n^T$ is equal to
the image of $D^T(K_n)$ in  $C^T(K_n)/C^T(K)$. We claim that $\phi_n^T$ is injective. For this
it clearly suffices to show the injectivity of the analogously defined map
$$ \phi_n : (\Z/\ell^n)^d \to  C(K_n)/C(K),$$
as $\phi_n$ factors through $\phi_n^T$. Let $(m_1,\ldots,m_d)
 \in \Z^d$ and assume that the class $[\prod_i b_i^{m_i}]$ in $C(K_n)$ is equal
 to $[\mathfrak c]$ where $\mathfrak c$ is a fractional ideal in $K$. This means
 that there exists $x \in ({K_n^\times}\otimes\Zell)^-$ such that ${\rm div}_{K_n}(x) = -\mathfrak c+\sum_i m_i\cdot b_i$.
Then the divisor on the right is $\Gamma_n$-invariant, hence $x^{\gamma_n-1}\in (O_{K_n}^\times\otimes\Zell)^-=\mu(K_n)\otimes\Zell$.
Since the module of roots of unity $\mu(K_n)$ is $\Gamma_n$-cohomologically trivial (well known fact),
we may
arrange that $x^{\gamma_n-1}=1$, that is $x$ is already in $K^\times\otimes\Zell$. Then the divisor $\sum_i m_i\cdot b_i$
also comes from $K$, and this is only possible if all $m_i$ are divisible by  $\ell^n$
(remember that all primes above $\ell$ in $K_n$ are totally ramified in $K_n/K$).
This shows that $\phi_n$ is injective as claimed. \ Q.E.D. \medskip

Recall that $C^T_\infty=\bigcup_n C^T(K_n)$.
Define $D^T_\infty:=\bigcup_n D^T(K_n)$.

\begin{lemma}\label{two}

(i)  $(C^T_\infty)^\Gamma = C^T(K) \cdot D^T_\infty$, and $D^T_\infty$ is divisible.

(ii)  $D^T_\infty$ is the divisible part of $(C^T_\infty)^\Gamma$.

(iii)  We have  $D^T_\infty \cap C^T(K) = D^T(K)$.

\end{lemma}

\Proof (i) We start with the ``ambiguous class number formula", both for $K_n/K$ and
for $K_n^+/K^+$, see Lemma 13.4.1 in  \cite{La}. If we divide
the former by the latter and note that the second factor in the denominator in loc.cit.
just goes away in the minus part (again, cohomological triviality of roots of unity),
we end up, after some comparison of notation, with the following:
$$  |C(K_n)^{\Gamma_n}| = |C(K)| \cdot \ell^{nd}. $$
It is a straightforward exercice to deduce from this the following $T$-variant:
$$  |C^T(K_n)^{\Gamma_n}| = |C^T(K)| \cdot \ell^{nd}. $$
When combined with the previous Lemma (and its proof) the above equality implies that the natural map $D_T(K_n)\to C^T(K_n)^{\Gamma_n}/C^T(K)$ is bijective.
Therefore we obtain
$$ C^T(K_n)^{\Gamma_n} = C^T(K) \cdot D^T(K_n). $$
By passing to the inductive limit, we obtain
$$  (C^T_\infty)^\Gamma = C^T(K) \cdot D^T_\infty. $$
This proves the equality in (i). Now, $D^T_\infty$ is divisible since all $\ell$-adic primes
are infinitely ramified in $K_\infty/K$.
Since $C^T(K)$ is finite, we get (ii) at once.

Part (iii) is proved using the method of proof of the preceding lemma: any element
of $D^T_\infty$ fixed by $\Gamma=\Gamma_0$ has to come from an ideal of $K$ supported
above $\ell$. \ Q.E.D.

\bigskip

We now present the initial step towards calculating
the coinvariants. We need one more object. Let $B^T(K)$
denote the quotient of $C^T(K)$ by the subgroup $D^T(K)$.

\begin{proposition}\label{old33ii}
There is an exact sequence
$$  0 \to B^T(K) \to \Tell(C^T_\infty)_\Gamma \to \Tell(\mathcal M)_\Gamma \to
\Zell[S\setminus S_\ell]^-\to 0. $$
\end{proposition}

\Proof We extract the following sequence from the diagram at the end of Section \ref{setup} (the second
arrow is the snake map):
$$  \Tell(\mathcal M)^\Gamma  \to \Zell[\cal S\setminus \cal S_\ell]^{-,\Gamma}
   \to \Tell(C^T_\infty)_\Gamma \to \Tell(\mathcal M)_\Gamma
   \to \Zell[\cal S\setminus\cal S_\ell]^-_\Gamma  \to 0. $$
The second and last nontrivial terms are isomorphic to  $\Zell[S\setminus S_\ell]^-$ (as primes in $S\setminus S_\ell$ are not ramified in $K_\infty/K$.) Going back to the proof of Proposition \ref{invariants},
one may verify the following: if we identify $\Tell(\mathcal M)^\Gamma $  with
$(\Zell \otimes U_{S,T})^-$ as in loc.cit, then the first arrow $\Tell(\mathcal M)^\Gamma  \to \Zell[S\setminus S_\ell]^{-}$
corresponds to the $S_\ell$--forgetful divisor map ${\rm div}_{K, S\setminus S_\ell}$ from
$(\Zell \otimes_{\Z} U_{S,T})^-$ to $\Zell[S\setminus S_\ell]^-$. Hence the cokernel
of the first arrow of the above sequence agrees with the cokernel of ${\rm div}_{K, S\setminus S_\ell}$;
this gives exactly $B^T(K)$, by definition, because of our assumption that $C^T(K)$ is
generated by $S$-ideal classes. The following commutative diagram, with surjective second row of vertical arrows
captures what is going on.
$$ \xymatrix{
     0 \ar[r]  & (U_{S_\ell, T}\otimes\Zell)^- \ar[r]\ar[d]  & (U_{S, T}\otimes\Zell)^-
     \ar[r]\ar[d]\ar[dr]^{{\rm div}_{K, S\setminus S_\ell}} & (U_{S, T}/U_{S_\ell, T}\otimes\Zell)^-  \ar[r]\ar[d] & 0  \\
0 \ar[r] & \Zell[S_\ell]^- \ar[r]\ar@{>>}[d]
      & \Zell[S]^- \ar[r]\ar@{>>}[d] &
     \Zell[S\setminus S_\ell]^- \ar[r]\ar@{>>}[d] & 0 \\
0 \ar[r] & D^T(K) \ar[r]  & C^T(K) \ar[r] &
     B^T(K)\ar[r]  & 0         }$$
This produces the exact sequence in the statement of part
(ii) of the lemma.
 \ Q.E.D. \medskip

Let $\alpha: B^T(K) \to \Tell(C^T_\infty)_\Gamma$ denote the first map
in the statement of the preceding proposition. We will determine the cokernel
of this map, and this will give the desired coinvariants. Let us state the result:

\begin{theorem}\label{coinvses} The following hold true.

(i) The cokernel of $\alpha: B^T(K) \to \Tell(C^T_\infty)_\Gamma$ is isomorphic to $\Zell[S_\ell]^-$.

(ii) We have a short exact sequence
$$  0 \to \Zell[S_\ell]^-  \to \Tell(\mathcal M)_\Gamma
   \to \Zell[S\setminus S_\ell]^-  \to 0. $$

(iii) We have a $\Zell[G]$--module isomorphism $T_\ell(\cal M)_{\Gamma}\cong \Zell[S]^-$.
\end{theorem}

\Proof Part (ii) is a direct consequence of Part (i) and Proposition \ref{old33ii}.
To prove part (i), we will need several lemmas. For simplicity, from this point on we will
let  $C:=C^T_\infty$. We will state the lemmas, explain why they suffice
to prove (ii) of the theorem, and then give the proofs of the lemmas. Then we will prove part (iii).

\begin{lemma}\label{new1}
There is an exact sequence
$$ 0 \to C^\Gamma/(C^\Gamma)_{div} \to \Tell(C)_\Gamma \to \Tell(C_\Gamma) \to 0. $$
\end{lemma}

\begin{lemma}\label{new2}
The left-hand term $C^\Gamma/(C^\Gamma)_{div}$ in Lemma \ref{new1} is isomorphic
to $B^T(K)$.
\end{lemma}

\begin{lemma}\label{new3}
The right-hand term $\Tell(C_\Gamma)$ in Lemma \ref{new1} is isomorphic to
$\Zell[S_\ell]^-$. In particular, it is torsion-free as a $\Zell$--module.
\end{lemma}

\medskip
Proof of Thm.~\ref{coinvses}(ii):  From Lemmas \ref{new1} and \ref{new3} we see that the torsion part of $\Tell(C)_\Gamma$
is exactly the image of the arrow $C^\Gamma/(C^\Gamma)_{div} \to \Tell(C)_\Gamma$.
Hence by Lemma \ref{new2} we infer that the torsion part of $\Tell(C)_\Gamma$
is isomorphic to $B^T(K)$. Now this is exactly the domain of definition of
the map $\alpha$. Even if we do not know the (injective) map $\alpha$,
we thus obtain that its cokernel identifies with the quotient of
$\Tell(C)_\Gamma$ modulo its torsion. Using the isomorphism of
Lemma \ref{new3}, we may conclude that the cokernel of $\alpha$
is isomorphic to $\Zell[S]^-$ as claimed. This concludes the proof of Thm.~\ref{coinvses}(i) and (ii) pending
the proofs of the lemmas. \ Q.E.D.

\bigskip

We now give the proof of the three lemmas in turn, the third one being
by far the most complex one. We tried to find a simpler argument, without
success. \medskip

\Proofof Lemma \ref{new1}:
Recall that our hypothesis that Iwasawa $\mu$--invariant associated to $K$ and $\ell$ vanishes implies
that $C$ is divisible. (See \cite{GP2} for details.) Therefore the short exact sequence of divisible groups
$$0 \to C/C^\Gamma \overset{1-\gamma}\longrightarrow  C \to
C_\Gamma \to 0$$ produces a short exact sequence of $\ell$--adic Tate modules
$$ 0 \to \Tell(C/C^\Gamma) \overset{1-\gamma}\longrightarrow \Tell(C) \to \Tell(C_\Gamma) \to 0. $$
Furthermore, noting that $\Tell(C)^\Gamma = \Tell(C^\Gamma)$, there is a canonical s.e.s.
$$0 \to \Tell(C)/\Tell(C)^\Gamma \to \Tell(C/C^\Gamma) \to C^\Gamma/(C^\Gamma)_{div} \to 0.$$
A diagram chase based on the two s.e.s.'s above then produces the desired s.e.s.
$$ 0 \to C^\Gamma/(C^\Gamma)_{div} \to \Tell(C)_\Gamma \to \Tell(C_\Gamma) \to 0. $$
Q.E.D. \medskip

\Proofof Lemma \ref{new2}:
  We need to calculate the quotient of $C^\Gamma$ by its maximal
	divisible subgroup. The latter is, by Lemma \ref{two} (ii),
equal to $D^T_\infty$. Hence $$C^\Gamma/(C^\Gamma)_{div}
= C^T(K) D^T_\infty/D^T_\infty \cong C^T(K)/C^T(K)\cap D^T_\infty
=C^T(K)/D^T(K) = B^T(K).$$ We used Lemma \ref{two} (i) and (iii).
\ Q.E.D.  \medskip

\Proofof Lemma \ref{new3}: We have to calculate the module $\Tell(C_\Gamma)$. As already mentioned,
this is the most delicate part.
We rely on Kurihara's paper \cite{Ku}, in particular on its Prop.~5.2,
which is proved using Lemma 5.1 of that paper. We apply this to the $\Gamma_n$-extension $K_n/K$,
and we note that we may omit the $\mu$-term at the left of the sequence in Prop.~5.2.
Kurihara's notation for the field extension is $L/K$; and we may omit the $\mu$-term
since it comes from a H$^1$ of
the $(-1)$-eigenspace of global units (first term in second
line of the long sequence in Lemma 5.1), so we may invoke cohomological
triviality of roots of unity again. Since all inertia groups of $K_n/K$ at
primes $v|\ell$ are the whole of $\Gamma_n$, the mentioned Proposition of [Ku]
gives the s.e.s.
$$ 0 \to \bigl( \bigoplus_{v\in S_\ell} \Gamma_n\bigr)^-
   \to  C(K_n)_{\Gamma_n} \to C(K)\to 0. $$
Routine arguments show that the following variant also holds:
$$ 0 \to \bigl( \bigoplus_{v\in S_\ell} \Gamma_n\bigr)^-
   \to  C^T(K_n)_{\Gamma_n} \to C^T(K) \to 0. $$
where the surjection is induced by the norm map at the level of ray class groups.
Since this norm map is onto, its kernel is isomorphic to
$\Hmone(\Gamma_n, C^T(K_n))$; on the other hand the term $\Gamma_n \cong \Htate(\Gamma_n,
K_{n,v}^\times)$, via the local Artin map. Consequently, we obtain an isomorphism
\begin{equation}\label{Kuri-iso} \Hmone(\Gamma_n, C^T(K_n))  \cong
\bigl( \bigoplus_{v\in S_\ell}
\Htate(\Gamma_n, K_{n,v}^\times)
\bigr)^- \end{equation} (see Kurihara's argument.) This isomorphism will be needed below.

We denote the norm map from $C^T(K_n)$ to $C^T(K)$ by $\pi_n$. Now, we need to pass to an inductive limit. To this end, we look at
the diagram
$$ \xymatrix{0 \ar[r] & \Hmone(\Gamma_n,C^T(K_n)) \ar[d]_{j_{n,n+1}} \ar[r]  & C^T(K_n)_\Gamma \ar[d] \ar[r]^{\pi_n} & C^T(K)
   \ar[d]_{\cdot \ell}  & \\
           0 \ar[r] & \Hmone(\Gamma_{n+1},C^T(K_{n+1}))  \ar[r]  & C^T(K_{n+1})_\Gamma  \ar[r]^{\pi_{n+1}} & C^T(K) .
}$$
Here the transition map $j_{n,n+1}$ has a direct and simple definition:
it is induced
by the inclusion map $C^T(K_n) \to C^T(K_{n+1})$ and the usual description
of $\Hmone$ as the kernel of the norm modulo the multiples of $(1-\sigma)$, with
$\sigma$ a generator of the cyclic group in question.
As $C^T(K)$ is finite, the inductive limit gives an isomorphism
$$  C_\Gamma \cong \varinjlim \Hmone(\Gamma_n,C^T(K_n)), $$
where the limit is taken along the maps $j_{n,n+1}$. Now \eqref{Kuri-iso} leads to an isomorphism
$$  C_\Gamma \cong \varinjlim \Htate(\Gamma_n, K_{n,v}^\times)$$
where the inductive
limit is taken along certain canonical maps
$$i_{n,n+1}: \Htate(\Gamma_n, K_{n,v}^\times)
\to \Htate(\Gamma_{n+1}, K_{n+1,v}^\times).$$
An easy direct
calculation reveals that $i_{n,n+1}$  is given by multiplication with the
relative norm element $\nu_{n+1,n}:=N_{G(K_{n+1}/K_n)}$. But in our case the
action of this element is the same as multiplication (or more properly,
exponentiation) by $\ell$. Therefore we
 have a commutative diagram:
$$\xymatrix{ \Htate(\Gamma_n, K_{n,v}^\times) \ar[r]  \ar[d]_{\nu_{n+1,n}=\ell} & \Gamma_n \ar[d]_{\ell}\\
 \Htate(\Gamma_{n+1}, K_{n+1,v}^\times) \ar[r]   & \Gamma_{n+1}
}$$
where the horizontal maps are local Artin maps. So we find that
$$C_\Gamma \cong \varinjlim \,(\bigoplus_{v\in S_\ell} \Gamma_n)^-,$$
where the transition maps are multiplication
by $\ell$. The choice of a generator for $\Gamma$ identifies the above injective limit
with $(\Q/\Z) \otimes_{\Z} \Zell[S_\ell]^-$. This proves,
by applying the functor $\Tell$, that $$\Tell(C_\Gamma) \cong
\Tell((\Q/\Z) \otimes \Zell[S_\ell]^-) \cong \Zell[S_\ell]^-,$$
which concludes the proof of Lemma \ref{new3}. \ Q.E.D.
\medskip

As mentioned earlier, the preceding series of arguments finishes the proof
of Theorem \ref{coinvses} parts (i) and (ii).
\bigskip

\Proofof of Thm.~\ref{coinvses} part (iii):
This would follow immediately if we could prove that the short exact sequence in  Thm.~\ref{coinvses}(ii)
is split. This is indeed the case, as shown
by the next lemma.

\begin{lemma}\label{perm-modules}
Let $G$ be any finite group and $U$ and $V$ two subgroups of $G$. Then
the Ext group $\Ext^1_{\Z [G]}(\Z[G/U],\Z[G/V])$ vanishes. (Our modules are
left modules, so $G/U$ denotes the set of left cosets $xU$, with the obvious
$G$-action.) Consequently $\Ext^1_{\Z[G]}(M,N)$ vanishes for any two
permutation modules $M$ and $N$, and this holds as well if the
base ring $\Z$ is replaced by $\Zell$.
\end{lemma}
\Proofof of Lemma \ref{perm-modules}: The proof is an exercise in permutation modules. We leave it
to the reader. \ Q.E.D.

The preceding Lemma applies in particular to the permutation
modules  $N=\Zell[S_\ell]^-$
and $M=\Zell[S\setminus S_\ell]^-$: the exact sequence in Thm.~\ref{coinvses} (ii) is split,
and the module in the middle is therefore isomorphic to $\Zell[S]^-$.
\ Q.E.D.

\begin{remark} Let us remark that we do not quite get
an explicit isomorphism between $\Tell(\mathcal M)_\Gamma$ and $\Zell[S]^-$. It is explicit up to a splitting
of an exact sequence, which exists but is not unique.  Unfortunately,
although the final theorem in
the next section does also imply, as a corollary,
that $\Tell(\mathcal M)_\Gamma$ is indeed isomorphic to $\Zell[S]^-$,
since the $\ell$-adified Tate sequence in the minus part has exactly $\Zell[S]^-$ on the
right, that isomorphism is much less explicit.\end{remark}


\section{Removing a technical assumption}\label{technical}

In this short section we explain how to eliminate condition (2) (see Section
\ref{Invariants}) in the end results of the preceding two sections.
(Condition (1) is built into the theory of Tate sequences and therefore indispensable.) The idea is the same
for invariants and for coinvariants. One chooses $n_0$ large enough
so that condition (2) holds for $K_\infty/K_{n_0}$ and
puts $\Gamma_0=\Gal(K_\infty/K_{n_0})$. If we replace $K$ by
$K_0$ in the results Prop.~\ref{invariants} and Thm.~\ref{coinvses} (ii)
(jointly with Lemma \ref{perm-modules}, we obtain descriptions
of $T_\ell(\cal M^{\Gamma_0})$ and $T_{\ell}(M_{\Gamma_0})$; the isomorphisms in
these descriptions are invariant under $G':=G \times (\Gamma/\Gamma_0)$.
We then perform a final (co)descent, taking
invariants (resp.~coinvariants) under the action $\Gamma/\Gamma_0$.
For the invariants everything is clear: the $\Gamma/\Gamma_0$-invariants
of $(U_{S,T}(K_{n_0})\otimes\Zell)^-$ coincide with $(U_{S,T}(K)\otimes\Zell)^-$. For the coinvariants,
it is also easy to check that $\Zell[S(K_{n_0})]^-_{\Gamma/\Gamma_0}$
is isomorphic to $\Zell[S(K)]^-$. The resulting isomorphism at level
$K$ is, of course, not quite explicit, since the isomorphism
at level $K_{n_0}$, coming about through Lemma \ref{perm-modules},
was not totally explicit.


\section{The link with Tate's canonical class}\label{tate-link}

We now consider the Tate canonical class $\tau := \tau_{K/k,S}\in \ $ Ext$^2_{\Z[G]}(X_S,U_S)$ introduced in Remark \ref{remark-ffields}.
We retain all our working hypotheses as well as notations introduced in Section 1.
In particular, $S$ is assumed large which means that the $S$--classes generate $cl_T(K)$, and consequently  $cl(K)$.
Tate proved (see \cite{Ta2}, Ch. 5, \S 2) that there exists a Yoneda $2$--extension  of $\Z[G]$-modules (not unique and not canonical)
\begin{equation}\label{Tate-sequence} \quad 0 \to U_S \to A \to B \to X_S \to 0\end{equation}
which represents $\tau$ and such that $A$ and $B$ are finitely generated and of finite projective dimension over $\Z[G]$ (i.e. cohomologically trivial or c.t. over $G$.)
Such a Yoneda extension is called a Tate sequence.
As mentioned before, we do not review the defining properties of $\tau$ here. The reader can consult
\cite{Ta1} and \cite{Ta2} for details.

It is our goal now to link $\tau$ with $\Tell(\cal M)$. For this, one
has to $\ell$-adify, $T$--modify and take the minus part of $\tau$, as explained in Section 1.
\medskip

Next, we follow \cite{BF} and \cite{BuDoc} and interpret the $\ell$--adification
$(\Zell\otimes_{\Z}\tau)$ of $\tau$ as the isomorphism class (in a sense to be made precise below) of the complex $[A_{\ell}\to B_{\ell}]$ in the derived category
$D^{\mbox{\tiny perf}}(\Zell[G])$ of perfect cochain complexes of $\Zell[G]$-modules,  where $A_\ell:=A\otimes\Zell$ and $B_\ell:=B\otimes\Zell$ are viewed
in degrees 0 and 1 respectively.
\medskip

Let $C^\bullet$ be a complex in the derived category $D(\Zell[G])$ (or $D(\Z[G])$) with differential maps $(\partial^s)_{s\in\Z}$ and some $i\in\Bbb Z$ such that
\begin{equation}\label{coh-vanishing} {\rm H^j}(C^\bullet)=0,\, \text{ for all } j\ne i, i+1.\end{equation}
Then one can associate to $C^\bullet$ the (correctly) truncated complex
$$\tau_{\geq i}(\tau_{\leq i+1} C^\bullet): [C^i/{\rm im\,}\partial^{i-1}\overset{\partial^i}\longrightarrow \ker\partial^{i+1}]$$
concentrated in degrees $i$ and $i+1$, with the same cohomology as $C^\bullet$.
This truncated complex leads to the canonical exact sequence
$$0\to{\rm H}^i(C^\bullet)\to (\tau_{\geq i}(\tau_{\leq i+1} C^\bullet))^i\to (\tau_{\geq i}(\tau_{\leq i+1} C^\bullet))^{i+1}\to {\rm H}^{i+1}(C^\bullet)\to 0,$$
which determines a Yoneda extension class $e(C^\bullet)\in {\rm Ext}^2_{\Zell[G]}({\rm H}^{i+1}(C^\bullet), {\rm H}^{i}(C^\bullet))$ (or ${\rm Ext}^2_{\Z[G]}$) canonically associated to
$C^\bullet$.

\begin{lemma}[Burns-Flach, \cite{BF}] Let $i\in\Bbb Z$ and $C^\bullet$ and $D^\bullet$ complexes in $D(\Zell[G])$ satisfying \eqref{coh-vanishing}.
Assume that we are given isomorphisms at the level of cohomology
$$\alpha_i: {\rm H^i}(C^\bullet)\overset{\sim}\longrightarrow {\rm H^i}(D^\bullet),\qquad \alpha_{i+1}: {\rm H^{i+1}}(C^\bullet)\overset{\sim}\longrightarrow {\rm H^{i+1}}(D^\bullet).$$
Then there exists an isomorphism $\alpha: C^\bullet\cong D^\bullet$ in $D(\Zell[G])$ such that $H^i(\alpha)=\alpha_i$ and $H^{i+1}(\alpha)=\alpha_{i+1}$ if and only if
$$(\alpha_{i+1}^{-1})^\ast\circ(\alpha_i)_\ast(e(C^\bullet))=e(D^\bullet),$$
where $(\alpha_{i+1}^{-1})^\ast\circ(\alpha_i)_\ast: {\rm Ext}^2_{\Zell[G]}({\rm H}^{i+1}(C^\bullet), {\rm H}^{i}(C^\bullet))\overset{\sim}\longrightarrow {\rm Ext}^2_{\Zell[G]}({\rm H}^{i+1}(D^\bullet), {\rm H}^{i}(D^\bullet))$
is the canonical isomorphism induced by $\alpha_i$ and $\alpha_{i+1}$.
\end{lemma}
\Proof See \cite{BF}, page 1353 or work out your own proof from the definitions. Q.E.D.

\begin{remark}
Note that for any complex $C^\bullet$ in $D(\Zell[G])$ satisfying \eqref{coh-vanishing} for some $i\in\Bbb Z$ there exists
an isomorphism in $D(\Zell[G])$
$$C^\bullet\cong \tau_{\geq i}(\tau_{\leq i+1} C^\bullet),$$
inducing the identity maps at the level of cohomology. So $e(C^\bullet)=e(\tau_{\geq i}(\tau_{\leq i+1} C^\bullet))$.

Most importantly, note that, by definition, any two Tate sequences
$$ 0 \to U_S \overset u\to A \overset f\to B \overset x\to X_S \to 0, \qquad  0 \to U_S \overset{u'}\to A' \overset {f'}\to B' \overset{x'}\to X_S \to 0 $$
give perfect complexes in $D^{perf}(\Zell[G])$ concentrated in levels $0$ and $1$
$$C^\bullet: [A_\ell\overset f\to B_\ell], \qquad C'^\bullet: [A'_\ell\overset {f'}\to B'_\ell]$$
and isomorphisms at the level of cohomology (induced by $u$, $u'$ and $x$, $x'$, respectively)
$${\rm H}^0(C^\bullet)\cong U_S\otimes \Zell\cong {\rm H}^0(C'^\bullet), \quad {\rm H}^1(C^\bullet)\cong U_S\otimes \Zell\cong {\rm H}^1(C'^\bullet)$$
which map the class $e(C^\bullet)$ to $e(C'^\bullet)$. Therefore, we have an isomorphism $C^\bullet\cong C'^\bullet$ in $D^{perf}(\Zell[G])$ which induces
the above isomorphisms at the level of cohomology.

\end{remark}

From now on we will denote by $(\tau\otimes\Zell)$ (respectively $(\tau\otimes\Zell)^-$) the complex $C^\bullet: [A_\ell\overset f\to B_\ell]$ (respectively $(C^\bullet)^-: [A_\ell^-\overset f\to B_\ell^-]$)
associated to a Tate sequence \eqref{Tate-sequence} as in the above remark. According to the above remark these complexes are unique up to isomorphisms in $D(\Zell[G])$.
\medskip

We will consider the affine schemes $$X:={\rm Spec}(O_{K})\setminus S={\rm Spec}(O_{K,S}),\qquad
\cal X:={\rm Spec}(O_{\cal K})\setminus\cal S={\rm Spec}(O_{\cal K,\cal S}).$$ We will
let $j:\cal T\to \cal X$ and $i:\cal X\setminus\cal T\to \cal X$ be the usual closed and open immersion,
respectively. When confusion is unlikely, we will use the same notation
$j: T\to X$ and $i: X\setminus T\to X$ for the corresponding immersions at the finite level. From now on
all cohomology is viewed in the \'etale sense, so in particular $R\Gamma(X, \ast):=R\Gamma(X_{et}, *)$,
$R\Gamma_c(X, \ast):=R\Gamma_c(X_{et}, *)$ and similarly for the scheme $\cal X$.

\begin{proposition}\label{tateclass-rgamma}  There is an isomorphism in $D^{\mbox{\tiny perf}}(\Zell[G])$
$$  (\Zell\otimes_{\Z} \tau)[-1] \cong R\Gamma(X, \Zell(1)). $$
(The $[-1]$--shift on the left produces a complex with cohomology
concentrated in degrees 1 and 2.)
\end{proposition}

\Proof  This is a fairly short argument.
All the same, it is not very direct, since it uses the
full strength of the key paper \cite{BF}.
Unexplained notation is taken literally from there;  all references in the present
proof are to this paper, if not said otherwise.

According to the last line of p.1383, the complex $\Psi_S$ represents
Tate's class $\tau$ (see the definition of $K_S$, p.1351 and p.1353).
By Prop.~3.3, we have an isomorphism
in the derived category (of course one also has to check, using the
explicit information given in loc.cit. eqn.(69) that
it gives the canonical maps on cohomology):
$$ \Zell\otimes \tau \cong R\Gamma_c(X,\Zell)^*[-2]. $$
The superscript star stands for
$R\,\Hom(-,\Zell)$ (a functor of the derived category to itself).
Now we invoke Lemma 16(b), which gives
$$ R\Gamma_c(X,\Zell)^* \cong  R\Gamma_c(X,\Q_\ell/\Zell)^\vee, $$
where the superscript $\vee$ is $R\Hom(-,\Q_\ell/\Zell)$. In contrast
to the functor $()^*$, the functor $\vee$ can be evaluated on {\it any\/}
complex in a quasi-isomorphism class, termwise, since $\Hom(-,\Q_\ell/\Zell)$
is exact. As a third and last ingredient, we invoke Artin-Verdier duality:
$$ R\Gamma_c(X,\Q_\ell/\Zell)^\vee[-3] \cong R\Gamma(X,\Zell(1)). $$
Again one has to make sure that the two preceding isomorphisms
are canonical on cohomology level.
Putting the three displayed isomorphisms together (the first shifted by
$-1$, and the second by $-3$), we obtain the formula of the proposition.
Note: We have been following the sign conventions of \cite{BF}. It appears
that in the terminology of \cite{BuDoc}, a minus sign would come up. \ Q. E. D.

\bigskip

Another important step is a description of $\Tell(\cal M)$ in terms of \'etale cohomology.
We intend to establish the following result.

\begin{theorem}\label{tatemotive-cohom} We have a canonical isomorphism
$$\Tell(\cal M) \cong \Hone(\cal X, j_{!}\Zell(1))^-. $$
\end{theorem}

\Proof We will actually prove $\cal M[m] \cong \Hone(\cal X, j_{!}\Z/m(1))^-$
for $m:=\ell^\nu$ and all $\nu\ge 1$. The isomorphisms will be compatible
and produce the desired result in the projective limit. We again resort to the
description given in \cite{GP2}:
$$\cal M[m] \cong
\bigl(\cal K^{(m)}_{\cal S,T}/\cal K_T^{\times m} \bigr)^-.$$
(See the proof of Proposition \ref{invariants} and the notations therein.)

\begin{proposition} For any fixed $m$ as above, there is a  natural isomorphism
$$ \phi=\phi_m: {\cal K}^{(m)}_{\cal S,T}/\cal K_T^{\times m} \overset\sim{\longrightarrow}
  \Hone(\cal X, j_{!}\Z/m(1)). $$
\end{proposition}

\Proof  This will take  several steps. Most of the underlying ideas
are from \cite{Del}, see Section 10.3.6 in particular, but the mathematical language in loc.cit. is
so different that we prefer to give a reasonably self-contained argument.

To make the main points more clearly visible, we will first prove
a simplified version: replace $T$ by the empty set. (In particular,
$j_{!}\Z/m(1)$ just becomes $\Z/m(1)$.) Then there is
an explicit geometric interpretation
of $\Hone(\cal X, \Z/m(1))$: it is canonically isomorphic to the group $D_m$
of equivalence classes of pairs $(\cal L,\alpha)$, where $\cal L$ is a projective rank one
module over $O_{\cal K,\cal S}$ (in other words a line bundle over $\cal X$),
and
$$   \alpha:  \cal L^{\otimes m}  \to  O_{\cal X} $$
is an isomorphism. The equivalence relation is as expected:
  $(\cal L,\alpha) \sim  (\cal L',\alpha')$  iff there is
an isomorphism $h: \cal L \to \cal L'$ with $\alpha'\circ h^{\otimes m}=\alpha$.
The group structure is obvious.
The relation between $D_m$ and $\Hone(\cal X, \Z/m(1))$ can be easily seen in the light
of Grothendieck's descent theory; the automorphism group
of the trivial element $(O_{\cal X},1)$ of $D_m$ is $\Z/m(1)$, just as the automorphism
group of the trivial (or any) line bundle is $\mathbb G_m$.

The isomorphism $\phi$ may now be constructed directly.
Given $f\in \cal K^{(m)}_{\cal S}$, we know that the principal $O_{\mathcal K}$-ideal
generated by $f$ is an $m$-th power away from $\cal S$, so
the sheaf $fO_{\cal X}$ is the $m$-th power of a unique ideal sheaf $\cal I$.
We let $\phi(\widehat f)$ be the class of the pair $(\cal I, f^{-1})$ in $D_m$. There are
two things to check: The kernel of
$\phi$ is precisely $\cal K^{\times m}$, and $\phi$ is surjective.
Both are straightforward. This settles the case where $T$ is replaced
by the empty set.

Now we put $T$ and $\cal T$ back in. (This is the part where our terminology and that in \cite{Del}
differ the most.) We define a modified group $D_m^{\cal T}$. Its elements are
equivalence classes of triples $(\cal L, \alpha, \beta)$, where $\cal L$ and $\alpha$ are as before and $\beta$ is defined as follows.
We let
$\kappa(\cal T):= \bigoplus_{v\in \cal T} \kappa(v),$
where $\kappa(v)$ is the residue field at $v$, as usual.
Now $\beta$, a so-called trivialization at $\cal T$, is
an isomorphism
$$\beta:  \kappa(\cal T) \otimes_{O_{\cal X}} \cal L  \overset\sim{\longrightarrow} \kappa(\cal T),$$
which has to be compatible with $\alpha$ in the obvious way:  $$id_{\kappa(\cal T)}\otimes \alpha =\beta^{\otimes m}.$$
Two triples $(\cal L, \alpha, \beta)$ and $(\cal L', \alpha', \beta')$ as above
are equivalent if there is an isomorphism $h: \cal L \overset\sim\to \cal L'$ such that
$$\alpha'\circ h^{\otimes m}=\alpha,\qquad \beta' = \beta \circ (id_{\kappa(\cal T)} \otimes_{O_{\cal X}} h).$$
The above argument carries over directly to
produce a canonical isomorphism between the groups ${\cal K}^{(m)}_{\cal S,T}/\cal K_T^{\times m}$ and $D_m^{\cal T}$.
It remains to identify $D_m^{\cal T}$ with \'etale cohomology. We feel this
should be known, and it certainly can be extracted from \cite{Del}
with some effort. Let us give a direct argument anyway, via
Cech cohomology.

Using that $m=\ell^\nu$ is invertible in $O_{\cal X}$, one easily obtains
that every element of $D_m^{\cal T}$ is trivialized by some \'etale covering
$(U_i)_i$ of $\cal X$. We may suppose that
all $U_i$ connected. The resulting transition maps over $U_i\cap U_j$
are on the one hand sections of $\Z/m(1)$ (as we said, this is
the automorphism sheaf of the trivial element of $D_m$), but because
of the trivialisation at $\cal T$ they are all trivial whenever $U_i\cap U_j$
has a point above $\cal T$. This produces therefore a 1-cocycle over the
sheaf $j_{!}\Z/m(1)$ relative to the covering, and hence a canonical map $\delta^{\cal T}$ from
$D_m^{\cal T}$ to the first Cech cohomology of that sheaf. Since Cech cohomology
embeds into \'etale cohomology, $\delta^{\cal T}$ gives a morphism $D^{\cal T}_m
\to \Hone(\cal X, j_{!}\Z/m(1))$. The analogous map with $\cal T$ empty is
an isomorphism. One has a commutative diagram
$$ \xymatrix{  1 \ar[r] &  \frac{\Z/m(1)(\kappa(\cal T))}{\Z/m(1)(\cal K)}   \ar[r]\ar[d]
        & D_m^{\cal T} \ar[r]\ar[d]^{\delta^T}  & D_m \ar[r] \ar[d]^{\delta} &  1 \\
 1 \ar[r]  & \frac{\Z/m(1)(\kappa(\cal T))}{\Z/m(1)(\cal K)}
       \ar[r]  &  \Hone(\cal X, j_{!}\Z/m(1)) \ar[r]  &  \Hone(\cal X, \Z/m(1)) \ar[r]  &  1. } $$
The top sequence comes from a standard s.e.s, cf. \cite{GPff}. One can
check directly that the leftmost vertical map is the identity.
Since $\delta$ is an isomorphism, $\delta^{\cal T}$ is an isomorphism as well. This proves the proposition. Q.E.D.

\medskip

Now, the proposition above together with the above mentioned identification of $\cal M[m]$ with
$\left ({\cal K}^{(m)}_{\cal S,T}/\cal K_T^{\times m}\right )^-$
and a passage to the projective limit, gives a proof of
Theorem \ref{tatemotive-cohom}. Q.E.D.

\medskip

In order to use the results \ref{tatemotive-cohom} and \ref{tateclass-rgamma}
towards our goal of identifying the Tate class
in terms of $\Tell(\cal M)$
we need some intermediate lemmas. All previous notation remains in place.

\begin{lemma}\label{onlyonedim} The  sheaf $j_{!}\Zell(1)$ on $\cal X$
has cohomology concentrated in degree 1.
\end{lemma}

\Proof To show this, one first looks at the cohomology of the $\ell$--adic \'etale sheaf sheaf $\Zell(1)$ on $\cal X$.

(1) $ \Hzero(\cal X, \Zell(1))=\varprojlim \mathbf{\mu}_{\ell^n}(\cal K)\cong \Zell(1)$ or $0$ if ${\mu}_\ell\subseteq\cal K^\times$ or not.

(2) $ \Htwo(\cal X, \Zell(1))=0$.  Indeed, if one writes the cohomology sequence attached to the Kummer sequence of \'etale sheaves on $\cal X$
$$0\to \Z/\ell^n(1)\to \Bbb G_m\to\Bbb G_m\to 0,$$
and takes into account that $\Hone(\cal X, \Bbb G_m)_\ell={\rm Pic}(O_{\cal K, \cal S})_\ell$ which is divisible under our working hypothesis that the $\mu$--invariant of $K$ and $\ell$ vanishes,
one concludes that $$\Htwo(\cal X, \Z/\ell^m(1))\cong{\rm Br}(\cal X)[\ell^m],$$
 for every $m$.
However, ${\rm Br}(\cal X)[\ell^m]=0$ for all $m$: if $A$ is a central simple algebra over $\cal K$ split outside $\cal S$ and killed
by $\ell^m$ (i.e. a representative of an element in ${\rm Br}(\cal X)[\ell^m]$) then it is defined over some $K_n$ and therefore split by $K_{n'}$, for $n'$ sufficiently large. (If $n'$ is sufficiently large,
the extension $K_{n'}/K_n$ has local degree divisible by $\ell^m$ at all primes in $S(K_n)\setminus S_\infty$ and therefore the algebra
$A$ is split by $K_{n'}$ locally everywhere and therefore splits globally.)  Passing to the limit gives the claimed vanishing.

Now we use the closed immersion $i: {\mathcal T} \to {\mathcal X}$ and the open immersion $j:\cal X\setminus \cal T\to \cal X$
and look at the standard exact sequence of sheaves on $\cal X$
$$  0 \to \js \Zell(1)  \to \Zell(1) \to  i_*\Zell(1) \to 0. $$
The long exact sequence in cohomology reads as follows.
$$ 0 \to \Hzero({\mathcal X} , \js\Zell(1)) \to \Hzero({\mathcal X}, \Zell(1))
\overset\rho{\longrightarrow} \Hzero({\mathcal X} , i_*\Zell(1))  \ \to \Hone({\mathcal X} , \js\Zell(1)) \to $$
$$  \Hone({\mathcal X}, \Zell(1))
\to \Hone({\mathcal X}, i_*\Zell(1)) \to \Htwo({\mathcal X}, \js\Zell(1))
\to \Htwo({\mathcal X}, \Zell(1))=0. $$
Now, the map $\rho$ is a diagonal embedding and therefore injective, as
$$\Hzero({\mathcal X} , i_*\Zell(1))\cong\oplus_{v\in\cal T}\Hzero(\kappa(v), \Zell(1))\cong \oplus_{v\in\cal T}\varprojlim \mu_{\ell^n}(\kappa(v))$$
and (under our working assumption on $T$) no roots of unity in $\mathcal K$ are congruent to $1$ mod $v$ for all $v\in\cal T$.
Consequently, $\Hzero({\mathcal X} , \js\Zell(1))=0$.
Now, $\cal T$ is a finite set of closed points on $\mathcal X$, so the natural map
$$\Hone({\mathcal X}, i_*\Zell(1)) \to
\Hone(\mathcal T, \Zell(1))$$
is an isomorphism. Since $\mathcal T$ is
a finite union of spectra of fields of char. $\ne \ell$ without algebraic extensions of $\ell$-power degree,
$\Hone(\mathcal T, \Zell(1))=0$. This implies (via the long exact sequence above) that $\Htwo({\mathcal X}, \js\Zell(1))=0$,
which concludes the proof. \ Q. E. D.

\medskip

For the purpose of the next results, we remind the reader that we are working under the hypothesis that $S$
is large (i.e. $cl_T(K)$ is generated by $S$--ideal classes.)

\begin{lemma}\label{tau-modif} \begin{enumerate}

\item The inclusion $\iota: U_{S,T} \to U_S$ induces a canonical isomorphism
$$\iota_\ast: \Ext^2_{\Z[G]}(X_S,U_{S, T})\cong \Ext^2_{\Z[G]}(X_S,U_{S}).$$

\item   The unique class $\tau':=\tau_{K/k, S, T}$ in $\Ext^2_{\Z[G]}(X_S,U_{S, T})$ satisfying  $\iota_\ast(\tau')=\tau$ admits a representative
$$0 \to U_{S,T} \to A' \to B '\to X_S \to 0$$
with $A'$ and $B'$ finitely generated and c.t. over $G$. (Any such representative will be called a $T$--modified Tate sequence.)

\item The pushout along $\iota$ of any $T$--modified Tate sequence is a Tate sequence.
\end{enumerate}
\end{lemma}

\Proof (1) Recall the exact sequence \eqref{basic-exact-seqeunce} in Section 1 and let $Z:=U_S/U_{S,T}$.
Since $S$ is large,  $Z\cong\kappa(T)$ as $\Bbb Z[G]$--modules. It is easily seen (see \cite{GP2})
that ${\rm pd}_{\Bbb Z[G]}(\kappa(T))=1$. Therefore $Z$ is c.t. over $G$. By a routine argument, we get that $\Ext^i_{\Z[G]}(N,Z)=0$ for all $i>0$,
all $G$-modules $N$ without $\Z$-torsion, and all $Z$ that are c.t. over $G$. This shows that the inclusion $\iota: U_{S,T} \to U_S$ induces an
isomorphism
$$\iota_*:{\rm Ext}^2_{\Z[G]}(X_S,U_{S,T})\cong {\rm Ext}^2_{\Z[G]}(X_S,U_S).$$

(2) Now, since $X_S$ is free of $\Z$--torsion, there is a canonical commutative diagram
$$ \xymatrix{
\Htwo(G, {\rm Hom}(X_S, U_{S,T}))\ar[d]^{\iota_\ast}_\wr\ar[r]^{\quad \sim} & {\rm Ext}^2_{\Z[G]}(X_S, U_{S,T})\ar[d]^{\iota_\ast}_\wr  \\
\Htwo(G, {\rm Hom}(X_S, U_{S}))\ar[r]^{\quad \sim} & {\rm Ext}^2_{\Z[G]}(X_S, U_{S})}
$$
Let $\alpha\in \Htwo(G, {\rm Hom}(X_S, U_{S}))$ be the preimage of $\tau$ via the bottom isomorphism. Then Tate showed (see \cite{Ta2}, Ch. II, \S5) that
the cup product with $\alpha$ induces isomorphisms $\widehat{\rm H}^i(G, X_S)\cong \widehat{\rm H}^{i+2}(G, U_S)$ , for all $i$.
Consequently, the cup product with $\alpha'$ induces similar isomorphisms $\widehat{\rm H}^i(G, X_S)\cong \widehat{\rm H}^{i+2}(G, U_{S, T})$. Now, this sufficient
for the argument in \cite{Ta2} pp. 56-57 (right before Remark 5.3 in loc.cit.) to produce a representative for $\tau'$ as required in part (2) of the Lemma.
One important note here is that since $U_{S,T}$ has no $\Z$--torsion (unlike $U_S$), $A'$ and $B'$ can be picked to be projective, finitely generated $\Z[G]$--modules.

(3) By definition, the push-out along $\iota$ of a $T$--modified Tate sequence as in (2) is a representative of $\tau$. It is of the
form $0\to U_S\to A\to B\to X_S\to 0$ with $B'=B$, hence f.g. and c.t. over $G$ and $A'$ part of an exact sequence
$0\to A'\to A\to \kappa(T)\to 0$, hence c.t. and f.g. over $G$. We obtain this way a Tate sequence.
\quad Q. E. D.
\medskip

\begin{lemma}\label{rgamma-modif}
  We have the following variant of Prop.~\ref{tateclass-rgamma}:
$$  (\Zell \otimes_{\Z} \tau')[-1] \cong R\Gamma(X,j_{!}\Zell(1)). $$
Here we have abusively used $j$ to indicate the open immersion $j: X\setminus T \to X$
at the finite level as well.
\end{lemma}

\Proof  Let  $\xi: \js \Zell(1)
\to \Zell(1)$ denote the canonical inclusion of sheaves.  Using the arguments in lemmas \ref{onlyonedim} and \ref{tau-modif} one checks
easily that H$^1(X,\xi)$ is injective with cokernel $U_S/U_{S,T} \cong \Hone(X,i^*\Zell(1))$,
and H$^2(X,\xi)$ is an isomorphism. It is then clear
that the $\ell$--adic \'etale sheaf $j_{!}\Zell(1)$ of $X$ has cohomology concentrated
in degrees 1 and 2 as well, so we can think of $R\Gamma(X, j_{!}\Zell(1))$
in terms of Yoneda 2-extensions.

Let $C^\bullet$ be a complex
concentrated in degrees 1 and 2 isomorphic in $D(\Zell[G])$ to $R\Gamma(X,j_{!}\Zell(1))$ . (Take for example the correct truncation of the latter complex.)
There is a map $f$ of complexes from $C^\bullet$  to some complex $D^\bullet$
which represents $R\Gamma(X,\Zell(1))$, such that $f$ induces H$^\bullet(X,\iota)$
on cohomology. In particular it gives the inclusion $U_{S,T} \to U_S$
on H$^1$, and an isomorphism on H$^2$. Let $\iota_*C^\bullet  ={C'}^\bullet$  be the complex given by pushing out:
$$ \xymatrix{
 0 \ar[r] &  \Zell \otimes U_{S,T} \ar[d]\ar[r] & C^1  \ar[d]\ar[r] & C^2  \ar@{=}[d]\\
 0 \ar[r] &  \Zell \otimes U_{S} \ar[r] & (C')^1 \ar[r]  & (C')^2.
   }$$
Note that ${C'}^2=C^2$. Then $f$ extends to a map of complexes $f'$ from ${C'}^\bullet$ to $D^\bullet$, just
by the universal property of the pushout. One verifies that $f'$ is now identity
on $\Hone$, and nothing has changed on $\Htwo$, so $f'$ is an isomorphism
and actually induces an equivalence. (See reminder before \ref{tateclass-rgamma}.)
Hence $\iota_*C^\bullet$ represents $R\Gamma(X,\Zell(1))$,
and this agrees with $\Zell\otimes \tau$ by Prop.~\ref{tateclass-rgamma}.
 Since $\tau'$ is
the inverse image of $\tau$ under $\iota_*$, we conclude
that $C^\bullet$ agrees with $\Zell \otimes \tau'$.\quad
Q.E.D.

\bigskip

With these preparations, we can state and prove the main result
of this section. Recall that $M:=\Tell(\cal M)$ and let $C^\bullet$ be the complex $[M \overset{1-\gamma}\longrightarrow  M]$
concentrated in degrees 0 and 1.

From this point on, for all $i\in\Z$ we let $M[i]$ denote the complex having $M$ in degree $(-i)$, $0$ everywhere else
and (obviously) $0$ differentials.

\begin{theorem}\label{mainresult} If the assumption at
the beginning of Lemma \ref{tau-modif} is satisfied,
then  there is a canonical isomorphism in $D(\Zell[G])$

$$  (\Zell \otimes_{\Z} \tau')^- \cong {C^\bullet}.$$

\end{theorem}

\Proof  By Theorem \ref{tatemotive-cohom} and Lemma \ref{onlyonedim} we have a canonical (therefore $\Gamma$--equivariant)
isomorphism in the derived category $D(\Zell[G])$
$$M[-1]\cong R\Gamma(\cal X, j_{!}\Zell(1))^-.$$


We
now descend from $\cal X$ to $X$. From \cite{BuDoc}, diagram (8) on p.371
plus comment (see also definition of $C(\theta)^\bullet$ on p.366 of loc.~cit.)
we get a canonical  isomorphism in the derived category $D(\Zell[G])$
\begin{equation}\label{etale-C}R\Gamma(X,j_{!}\Zell(1))^-\cong C^\bullet[-1].\end{equation}
Three comments are necessary in order to derive this isomorphism from loc. cit.

(1) To link up with the
notation in \cite{BuDoc}, note that the $(-1)$ shift of the mapping cone
of the map of complexes $1-\gamma: M[-1] \to M[-1]$ (which is the precise definition of Burns' $C(\theta)^\bullet$ in our context)
is exactly the complex $C^\bullet[-1]$.

(2) We also remark that
\cite{BuDoc} is concerned with the function field
case where there is a canonical choice for $\gamma$, to wit Frobenius.
But actually the isomorphism class (in the derived sense) of the complex $C^\bullet: [M\stackrel{1-\gamma}{\rightarrow} M]$
does not change
when $\gamma$ is replaced by any other generator of $\Gamma$, so the
lack of a canonical generator of $\Gamma$ is not an issue.

(3) The rest of the argument taken from \cite{BuDoc} is entirely
cohomological algebra, so there is no difference between
the function field and number field cases in this respect.

Finally, we combine \eqref{etale-C} with Lemma \ref{rgamma-modif} to get the isomorphism
in the statement of the above theorem.\quad  Q. E. D.

\medskip

\begin{remark}\label{functionfields} (1) If one weakens the assumption
at the beginning of Lemma \ref{tau-modif} to say that just
$C^T(K)$ (the minus-$\ell$-part of $cl(K)$) is $S$-generated (which
is even closer to Hypothesis (1) in Section \ref{Invariants}), then
an $\ell$-adic version of that lemma remains correct, as well
as a version of the preceding theorem, in which the modified
Tate class $\tau'$ only exists as an $\ell$-adic object.

(2) With notations as in the proof of Lemma \ref{rgamma-modif}, one also has
$(\Zell \otimes_{\Z} \tau)^- \cong ({C'}^\bullet)^-$. This can be seen as an ``explicit Tate sequence".
The complex ${C'}^\bullet$ arises from $M$ by a very simple and
explicit construction involving pushout along $\iota$.
For further reference, here is the relevant diagram (basically
taken from the proof of Lemma \ref{rgamma-modif}; we also put
in the cokernels on the right for clarity);
$M'$ is defined as the pushout,
and $({C'}^{\bullet})^-$ is simply the complex $[M'\to M]$ that
shows up in the lower row.
$$ \xymatrix{
 0 \ar[r] &  (\Zell \otimes U_{S,T})^- \ar[d]\ar[r] & M  \ar[d]\ar[r]^{1-\gamma} & M  \ar@{=}[d]\ar[r]
        & \Zell[S]^- \ar[r]\ar@{=}[d] & 0 \\
 0 \ar[r] &  (\Zell \otimes U_{S})^- \ar[r] & M' \ar[r]  & M  \ar[r] & \Zell[S]^- \ar[r] & 0.
   }$$
Note that in order to really work with $[M'\to M]$, one needs a good grasp
on the maps $(\Zell \otimes U_{S,T})^- \to M$ and $M \to \Zell[S]^-$. This is
another  justification, apart from their intrinsic interest,
for the explicit calculations in Sections \ref{Invariants} and \ref{Coinvariants}.
\end{remark}

\begin{remark}\label{linear-disjointness-removed}
Finally, we would like to indicate briefly how the linear disjointness condition $k_\infty\cap K=k$ can be removed
in all of the above considerations. In the case where this condition is not satisfied,
$T_\ell(\mathcal M_{S,T}(K_\infty))$ does not have a natural $\Z_{\ell}[G]$--module structure. Indeed, in this case
$G(K_\infty/k)\cong H\rtimes\Gamma$, where $H:=G(K/k_\infty\cap K)$ and $\Gamma:=G(k_\infty/k)$, so $T_\ell(\mathcal M_{S,T}(K_\infty))$
is naturally endowed with a $\Bbb Z_\ell[H]$--module structure only and it is projective over this ring (see \cite{GP2}.). Consequently
$T_\ell(\mathcal M_{S,T}(K_\infty))\otimes_{\Zell[H]}\Zell[G]$ is a projective $\Zell[G]$--module. It is easily seen (see \cite{GP2}) that this is in fact isomorphic to
the $\ell$--adic realization of the abstract $\ell$--adic $1$--motive associated to the semisimple $k$--algebra $K\otimes_k k_\infty$ and the sets $S$ and $T$, i.e.
we have a natural isomorphism of $Z_\ell[[G\rtimes\Gamma]]$--modules
$$T_\ell(\mathcal M_{S,T}(K\otimes_k k_\infty))\cong T_\ell(\mathcal M_{S,T}(K_\infty))\otimes_{\Zell[H]}\Zell[G].$$
All of the above considerations can be easily generalized to show that the complex
$$[T_\ell(\mathcal M_{S,T}(K\otimes_k k_\infty))\overset{1-\gamma}\longrightarrow T_\ell(\mathcal M_{S,T}(K\otimes_k k_\infty))]$$
concentrated at levels $0$ and $1$ represents the (minus $\ell$--adic) Tate class and gives an explicit (minus $\ell$--adic) Tate sequence
just as above.
\end{remark}


\section{Examples}

We finish this paper by sketching one or two examples, without going into
detail too deeply. The main purpose is twofold: to have a certain reality check
on our results, and to give the reader a feeling what is going on.

We choose a setting that is as simple as possible. Let $k=\Q$ and
$K^+$ the cubic field of conductor 7. We take $\ell=3$. For $K$, we
will look at two choices: $K=K^+L$ where the
imaginary quadratic field $L$ is either $\Q(\sqrt{-5}$
or $\Q(\sqrt{-37})$.
In both cases $K$ is CM and of course $k$ is totally real. For $S$ we
consistently take the set of ramified primes in $K/k$ together with
the 3-adic primes; the contribution
of the places over $5$ (resp.~37) and of the infinite places
disappears in the minus part. For
$T$ we take the set of places in $K$ above any totally split place
in $K/k$. In both cases, 7 is split in $L$
and ramified in $K/L$. Moreover $\zeta_3$ is not
contained in $K_\infty$. Hence the ``toric part", that is,
the kernel of $T_\ell(C^T_\infty) \to T_\ell(C_\infty)$, is
a copy of the free module $\Zell[G']$, where $G'=\Gal(K/L)=\Gal(K^+/\Q)$.
The rank of the lattice $\Zell[S-S_\ell]^-$ is 1 in both cases.
Thus we have in both cases, recalling that $M=T_\ell(\mathcal M)$:
$$ rk(M) = \lambda_{3,K}^- + 4, $$
where the constant 4 comes about as $3+1$; 3 for the toric part and
1 for the lattice part. Since $M$ is free over $\Zell[G']$, this
already tells us that $\lambda_{3,K}^- \equiv 2$ modulo 3.

This can also be seen from the Kida formula which says
$$  \lambda_{3,K}^- = 3 \lambda_{3,L} + 2. $$

First case: $K=K^+(\sqrt{-5})$. Here 3 is split in $L$,
and by \cite{DFKS}, $\lambda_{3,L}=
\lambda_{3,L}^-=1$. Hence $\lambda_{3,K}^-=5$ and the
rank of $M$ is 9. Both the $\Gamma$-invariants and coinvariants
of $M$ give a rank 2 module with trivial $G'$-action.

Second case: $K=K^+(\sqrt{-37})$. Here 3 is inert in $L$,
and therefore $\lambda_{3,L}=\lambda_{3,L}^-=0$.
Hence $\lambda_{3,K}^-=2$ and the
rank of $M$ is 6. Both the $\Gamma$-invariants and coinvariants
of $M$ give a rank 1 module with trivial $G'$-action.

Final remark on the first case: Since the toric part has no $\Gamma$-invariants,
and the lattice part of $M$ has rank one, it follows that the
module of $\Gamma$-invariants in $M_0:=T_\ell(C_\infty)$ has rank one.
This already excludes that the rank of $M_0$ is 2,
as happens in the second case. Indeed, if the rank were
2, then $M_0$ would be annihilated by the norm element of $G'$,
so $M_0$ would be free over the DVR $\Z_3[G']/N_{G'}$. We know that
a chosen generator
$\gamma$ of $\Gamma$ has an eigenvalue 1 on this module; but then the
characteristic polynomial of $\gamma$ would have to be $(x-1)^2$
(in other words, the eigenvalue 1 would have algebraic multiplicity 2).
This would contradict the 3-adic Gross conjecture, which
states that the quotient of $T_\ell(C_\infty)$ by its $\Gamma$-invariants
has no $\Gamma$-invariants.

\medskip

\noindent{\bf Acknowledgment:} The first author would like to thank
Guido Kings for some very helpful hints.


\end{document}